\documentclass{article}
\usepackage{amssymb,amsmath,mathrsfs,latexsym,amscd}
\usepackage{exscale,latexsym,amsthm,graphics}

\usepackage{makeidx}
\makeindex
\usepackage{epsfig}
\usepackage{txfonts}
\newtheorem{thm}{Theorem}[section]

\newtheorem{lemma}[thm]{Lemma}
\newtheorem{prop}[thm]{Proposition}

\newtheorem{defn}[thm]{Definition}

\newtheorem{remark}[thm]{Remark}

\def\bM {{\mathbb M}}  
 \def\Q {{\mathbb Q}}

\def\N {{\mathbb N}}

\def\qed{{\hfill $\Box$ \bigskip}}
\def\N {{\mathbb N}}
\def\R {{\mathbb R}}

\def\P{{\mathbb P}}

\def\E{{\mathcal E}}
\def\P{{\mathbb P}}

\begin{document}

\noindent
{{\Large\bf Strict decomposition of diffusions associated to degenerate (sub)-elliptic forms }}

\bigskip
\noindent
{\bf Jiyong Shin}\\

\noindent
{\small{\bf Abstract.} For a given strongly local Dirichlet form with possibly degenerate  symmetric (sub)-elliptic matrix, we identify a Hunt process (associated to the Dirichlet form) with a weak solution to  the  corresponding stochastic differential equation  starting from all points in $\R^d$. More precisely, using heat kernel estimates, stochastic calculus, and Dirichlet form theory, we obtain the pointwise existence of a weak solution to the stochastic differential equation which has possibly unbounded and discontinuous drift. We also present some conditions that the weak solution becomes a pathwise unique strong solution.\\ 



\section{Introduction}
In this paper, we are concerned with a symmetric Dirichlet form (given as the closure of)
\begin{equation}\label{subdf}
\E^D(f,g): = \frac{1}{2} \int_{\R^d} \langle D \ \nabla f , \nabla g \rangle \ dx, \quad f,g \in C_0^{\infty}(\R^d)
\end{equation}
on $L^2(\R^d,dx)$ and the corresponding stochastic differential equation (hereafter SDE)
\begin{equation}\label{sdedd}
X_t = x + \int_0^t \sigma(X_s) \ dW_s + \int_0^t b(X_s) \ ds, \quad x \in \R^d,
\end{equation}
where the conditions on the (possibly) degenerate diffusion matrix $D = (d_{ij})_{1 \le i,j \le d}$ are formulated in (A1) and (A2) of Section \ref{sepo2} and in (A3), (A4), and (A5)$^{\prime}$ of Section \ref{se3deel} (for $\sigma$, $b$ see Theorem \ref{th;exsolsu}  and Theorem \ref{t;sdpoe}).

Given the bilinear form \eqref{subdf}, it is well known from Dirichlet form theory (Fukushima decomposition) and localization method  that one may derive a weak solution to the SDE \eqref{sdedd} for any starting point $x \in \R^d \setminus N$, where $N$ is some capacity zero set w.r.t. $\E^D$ (see \cite{FOT}). For the Dirichlet form $\E^D$ with uniformly elliptic matrix D, it has been shown  in \cite[Example]{fuku93}) that the weak solution to the corresponding SDE \eqref{sdedd} exists for all starting points in $\R^d$. However, in general for the Dirichlet form $\E^D$ with the possibly degenerate matrix $D$ there is  no characterization of $D$ which allows to give rise to a weak solution to the corresponding SDE \eqref{sdedd} for explicitly specified starting points in $\R^d$.

In this point of view, the main aim of this article is to construct a Hunt process associated to $\E^D$ (degenerate (sub)-elliptic form) which satisfies the Fukushima's absolute continuity condition (cf. \cite[(4.2.9) and Theorem 5.5.5]{FOT}) and in the sequel to identify it with the solution of the associated SDE \eqref{sdedd} for any starting point $x \in \R^d$.  The identification of the process with the solution to the SDE \eqref{sdedd} starting from explicitly specified points is of central interest in Dirichlet form theory. Following the tools and techniques developed in  \cite{ShTr13a} and \cite{ShTr15} we construct a Hunt process satisfying the absolute continuity condition and identify it with the solution to the SDE \eqref{sdedd} pointwise under some additional assumptions, namely (A1), (A2) in Section \ref{sepo2} and (A3), (A4), and (A5)' in Section \ref{se3deel}.  In \cite{ShTr13a} and \cite{ShTr15}, the (strong) equivalence between the intrinsic metric (derived from the Dirichlet form $\E$ there) and the Euclidean metric plays a crucial role throughout the articles.  In this paper we show that  the  local equivalence between the intrinsic metric (derived from $\E^D$) and the Euclidean metric is enough to obtain similar results (see \eqref{eq;intldis}, \eqref{eq;intudisloc} and Lemma \ref{le;l;intrinsicdell}). Therefore this paper is basically a continuation of \cite{ShTr15}. To our knowledge, however, it is first time to consider the essential degenerate matrix $D$ in the bilinear form \eqref{subdf} and via Dirichlet form theory to show the existence of a weak solution to the corresponding SDE for any starting points in $\R^d$ (see Section \ref{sepo2}).
       
The contents of this paper are organized as follows. In Section \ref{sepo2}, we consider a symmetric diffusion matrix  $A=(a_{ij})_{1 \le i,j \le d}$ satisfying the subelliptic estimate. We first present analytic background based on the results from \cite{BM, FePh, Je, JeSa, NSW, St3, St5}. Then using local equivalence between the intrinsic metric and the Euclidean metric, we show that the Dirichlet form $(\E^A, D(\E^A))$ is conservative and even recurrent in the case of $d=2$ (see Theorem \ref{t;subconsv}). 
In order to construct the Hunt process associated with $(\E^A, D(\E^A))$ satisfying the absolute continuity condition we apply the Dirichlet form method developed in \cite{ShTr13a} and finally identify it with the solution to  the SDE \eqref{sdedd}. In Section \ref{se3deel} we consider a different degenerate (locally uniformly) elliptic matrix $B$
and do the same as in Section \ref{sepo2}. In this case, however, unlike Section \ref{sepo2} we can show that the associated semigroup is Feller in classical sense.  Section \ref{sect4pust} is devoted to pathwise uniqueness and strong solution.\\

\noindent{\bf{Notations}}:\\
For an open set $E \subset \R^d$, $d \ge 2$ with Borel  $\sigma$-algebra $\mathcal{B}(E)$ we denote the set of all $\mathcal{B}(E)$-measurable $f : E \rightarrow \R$ which are bounded by $\mathcal{B}_b(E)$. The usual $L^q$-spaces $L^q(E, \mu)$, $q \in[1,\infty]$ are equipped with $L^{q}$-norm $\| \cdot \|_{L^q(E,\mu)}$ with respect to the  measure $\mu$ on $E$, $\mathcal{A}_b$ : = $\mathcal{A} \cap \mathcal{B}_b(E)$ for $\mathcal{A} \subset L^q(E,\mu)$,  and $L^{q}_{loc}(E,\mu) := \{ f \,|\; f \cdot 1_{U} \in L^q(E, \mu),\,\forall U \subset E, U \text{ relatively compact open} \}$, where $1_A$ denotes the indicator function of a set $A$. If $\mathcal{A}$ is a set of functions $f : E \to \R$, we define $\mathcal{A}_0 : = \{f \in \mathcal{A} \ | $ supp($f$) : = supp($|f| \mu$) is compact in $E \}$. Let $\nabla f : = ( \partial_{1} f, \dots , \partial_{d} f )$ where $\partial_j f$ is the $j$-th weak partial derivative of $f$ and $\partial_{ij} f := \partial_{i}(\partial_{j} f) $, $i,j=1, \dots, d$.  We denote the set of continuous functions on $E$, the set of continuous bounded functions on $E$, the set of compactly supported continuous functions in $E$ by $C(E)$, $C_b(E)$, $C_0(E)$, respectively.  The space of continuous functions on $E$ which vanish at infinity is denoted by $C_{\infty}(E)$. The set of all infinitely differentiable functions on $E$ and the set of all infinitely differentiable functions with compact support in $E$ are denoted by $C^{\infty}(E)$ and $C_0^{\infty}(E)$, respectively.  As usual we denote the Lebesgue measure on $\R^d$ by $dx$ and equip $\R^d$ with the Euclidean norm $\| \cdot \|$ and the corresponding inner product $\langle \cdot, \cdot \rangle$. We write $B_{r}(x): = \{ y \in \R^d \ | \ \|x-y\| < r  \}$, $x \in \R^d$, $r>0$. For a Borel measurable set $A \subset \R^d$, the closure of $A$ in $\R^d$ is denoted by $\overline{A}$ and the volume of A w.r.t. Lebesgue measure is denoted by $|A|$.

\section{Preliminaries and degenerate subelliptic forms with Lebesgue measure}\label{sepo2}

Throughout this paper, we consider a symmetric matrix $D=(d_{ij})_{1 \le i,j \le d}$, $d_{ij} \in L^1_{loc}(\R^d, dx)$ and a symmetric bilinear form
\[
\E^{D}(f,g) : = \frac{1}{2} \int_{\R^d} \langle D \ \nabla f , \nabla g \rangle \ dx, \quad f,g \in C_0^{\infty}(\R^d).
\]
Suppose to define some notations that the symmetric bilinear form $(\E^{D},C_0^{\infty}(\R^d))$ is closable in $L^2(\R^d, dx)$ and its closure $(\E^D,D(\E^D))$ is a strongly local, regular, symmetric Dirichlet form (cf. \cite{FOT}). The Dirichlet form $(\E^D,D(\E^D))$ can be written as 
\[
\E^D(f,g) = \frac{1}{2} \int_{\R^d} \ d\Gamma^{D} (f,g), \quad f,g \in D(\E^D),
\]
where $\Gamma^{D}$ is a symmetric bilinear form on $D(\E^D) \times D(\E^D)$ with values in the signed Radon measures on $\R^d$ (called energy measures).
The nonnegative definite measure $\Gamma^D(f,f)$ can be defined by the formula
\[
\int_{\R^d} \phi \ d\Gamma^D(f,f) = 2 \E^D (f,\phi f) -  \E^D(f^2,\phi),
\]
for every $f \in D(\E^D) \cap L^{\infty}(\R^d,dx)$ and every $\phi \in D(\E^D) \cap C_0(\R^d)$. Let $D(\E^D)_{loc}$ be the set of all measurable functions $f$ on $\R^d$ for which on every compact set $K \subset \R^d$ there exists a function $g \in D(\E^D)$ with $f=g$ $dx$-a.e on $K$. By an approximation argument we can extend the quadratic form $f \mapsto \Gamma^{D}(f,f)$ to $D(\E^D)_{loc} = \big\{ f \in L^2_{loc}(\R^d, dx) \ | \ \Gamma^{D}(f,f) \ \text{is a Radon}\big.$ $\big. \text{measure} \big\}$. By polarization we then obtain for $f,g \in D(\E^D)_{loc}$ a signed Radon measure
\[
\Gamma^D(f,g) = \frac{1}{4} \Big( \Gamma^D(f+g,f+g)) - \Gamma^D(f-g,f-g) \Big).
\]
For these properties of energy mesures we refer to \cite{FOT}, \cite[Proposition 1.4.1]{Le}, and \cite{Mo} (cf. \cite[Appendix]{St3}). The energy measure $\Gamma^{D}$ defines in an intrinsic way a pseudo metric $d$ on $\R^d$ by
\begin{equation}\label{inmetc}
d(x,y) = \sup \Big\{f(x) - f(y) \ | \ f  \in D(\E^D)_{loc} \cap C(\R^d), \ \Gamma^{D}(f,f)   \le   dz \ \text{on} \ \R^d  \Big\}, 
\end{equation}
where $\Gamma^{D}(f,f)   \le   dz$ means that the energy measure $\Gamma^{D}(f,f)$ is absolutely continuous w.r.t. the reference measure $dz$ with Radon-Nikodym derivative $\frac{d \Gamma^{D}(f,f)}{dz}  \le 1$ (cf. \cite{BM}).  We define the balls w.r.t. the intrinsic metric by 
\[
\tilde{B}_r(x) = \{ y \in \R^d \ | \ d(x,y) < r \}, \quad x \in \R^d, \quad r>0.
\]
\begin{defn}
\begin{itemize}
\item[(i)]
We say the completeness property holds, if for all balls $\tilde{B}_{2r}(x)$,  $x \in \R^d$, $r>0$, the closed balls $\overline{\tilde{B}_r}(x)$ are complete (or equivalently, compact).
\item[(ii)]
We say the doubling property holds for a given measure $\mu$, if there exists a constant N= N(d) such that for all balls $\tilde{B}_{2r}(x) \subset \R^d$
\[
\mu (\tilde{B}_{2r}(x)) \le 2^N \mu (\tilde{B}_{r}(x)).
\] 
\item[(iii)] We say the (scaled) weak Poincar\'e inequality holds, if there exists a constant $c= c(d)$ such that for all balls $\tilde{B}_{2r}(x) \subset \R^d$
\[
\int_{\tilde{B}_r(x)} |f - \tilde{f}_{x,r} |^2 \ dy \le c \ r^2 \int_{\tilde{B}_{2r}(x)} \ d\Gamma^{D}(f,f), \quad f \in D(\E^D),
\]
where  $\tilde{f}_{x,r} = \frac{1}{|\tilde{B}_r(x)| }\int_{\tilde{B}_r(x)} f \ dy$.
\item[(iv)] A strongly local, symmetric Dirichlet form $(\E^D,D(\E^D))$ is called strongly regular if it is regular and if $d(\cdot,\cdot)$ (defined by \eqref{inmetc}) is a metric on $\R^d$ whose topology coincides with the original one.
\end{itemize}
\end{defn}

A positive Radon measure $\mu$ on $\R^d$ is said to be of finite energy integral if
\[
\int_{\R^d} |f(x)|\, \mu (dx) \leq c \sqrt{\E^{D}_1(f,f)}, \; f\in D(\E^{D}) \cap C_0(\R^d),
\]
where $c$ is some constant independent of $f$ and $\E^{D}_1(f,f) : = \E^{D}(f,f) + \int_{\R^d} |f|^2  \ dx$. A positive Radon measure $\mu$ on $\R^d$ is of finite energy integral, if and only if there exists a unique function $U_{1} \, \mu\in D(\E^{D} )$ such that
\[
\E^{D}_{1}(U_{1} \, \mu, f) = \int_{\R^d} f(x) \, \mu(dx),
\]
for all $f \in D(\E^{D}) \cap C_0(\R^d)$. $U_{1} \, \mu$ is called the $1$-potential of $\mu$. The measures of finite energy integral are denoted by $S_0$. We further define $S_{00}: = \{\mu\in S_0 \, | \; \mu(\R^d)<\infty, \|U_{1} \mu\|_{ L^{\infty}(\R^d,dx)} <\infty \}$.\\

In this section, we consider the following assumption:
\begin{itemize}
\item[(A1)] $A=(a_{ij})_{1 \le i,j \le d}$ is a symmetric matrix such that 
\[
a_{ij} \in C^{\infty} (\R^d) \cap C_b(\R^d), \quad i,j=1, \dots, d,
\]
and A satisfies the degenerate elliptic condition (positive semidefinite), i.e. for  $dx$-a.e. $x \in \R^d$
\begin{equation*}
0 \le \langle A(x) \ \xi, \xi \rangle,  \quad \forall \xi \in \R^d.
\end{equation*}
\end{itemize}
From now on we fix a symmetric matrix $A=(a_{ij})_{1 \le i,j \le d}$ satisfying (A1) and consider the symmetric bilinear form
\[
\E^{A}(f,g) : = \frac{1}{2} \int_{\R^d} \langle A \ \nabla f , \nabla g \rangle \ dx, \quad f,g \in C_0^{\infty}(\R^d).
\]
Furthermore we assume:
\begin{itemize}
\item[(A2)]
The symmetric matrix $A$ satisfies the following subelliptic estimate, i.e. there exist constants $\varepsilon > 0$ and  $\delta >0$
such that 
\begin{equation*}
\delta \ \| u \|^2_{H^{\varepsilon}} \le \E^{A} (u,u) + \|u\|^2_{L^2(\R^d,dx)}, \quad \forall u \in C_0^{\infty}(\R^d).
\end{equation*}
Here $\| u \|^2_{H^{\varepsilon}}:= \int_{\R^d} |\hat{u} (\xi)|^2 \cdot (1+ \|\xi \|^2 )^{\varepsilon} \ d \xi$ for any $\varepsilon > 0$ and $\hat{u}$ is the Fourier transform of $u$ and $H^{\varepsilon}(\R^d):= \{ u \in L^2(\R^d, dx) \ | \ \|u\|_{H^{\varepsilon}} < \infty \}$ is the fractional Sobolev space of order $\varepsilon > 0$. 
\end{itemize}
\begin{remark}
We refer to \cite[(1.3), Theorem 2.1, 2.2]{JeSa} for some operators satisfying the subelliptic estimate (A2).
\end{remark}
 By \cite[Section 3.1 (1$^{\circ}$)]{FOT} $(\E^{A},C_0^{\infty}(\R^d))$ is then closable in $L^2(\R^d, dx)$ and its closure $(\E^A,D(\E^A))$ is a strongly local, regular, symmetric Dirichlet form. Furthermore, it is known that that the intrinsic metric $d(\cdot,\cdot)$ derived from the Dirichlet form $(\E^A,D(\E^A))$ satisfies (see \cite[Theorem 4.2]{St5})
\begin{equation}\label{eq;intldis}
d(x,y) \ge c_0^{-1} \|x-y\|, \quad \forall x,y \in \R^d,
\end{equation}
where $c_0 \ge 1$ is some constant (see also the proof of Lemma \ref{le;l;intrinsicdell} below) and there exist $r_0 > 0$, $C_0 >0$ such that
\begin{equation}\label{eq;intudisloc}
d(x,y)  \le C_0 \|x-y\|^{\varepsilon}, \quad \forall x,y \in \R^d \ \text{with} \ \|x-y\| < r_0,
\end{equation} 
where $\varepsilon \in (0,1)$ is the constant as in (A2) (see  \cite[Section 1. (b)]{BM} and \cite{FePh, JeSa}). Let  $(T_t)_{t > 0}$ and $(G_{\alpha})_{\alpha > 0}$ be the $L^2(\R^d, dx)$-semigroup and resolvent associated to $(\E^A,D(\E^A))$ (see \cite{FOT}).

\begin{remark}
The topology induced by the intrinsic metric coincides with the Euclidean topology (see  \cite[Section 1. (b)]{BM} and \cite{NSW}). 
Hence the Dirichlet form $(\E^A,D(\E^A))$ is strongly regular. In particular $\tilde{B}_r(x) \in \mathcal{B}(\R^d)$ for all $x \in \R^d, r>0$.
\end{remark}
\begin{thm}\label{t;subconsv}
\begin{itemize}
\item[(i)] Let $d=2$. Then the Dirichlet form $(\E^A,D(\E^A))$ is recurrent.
\item[(ii)] The Dirichlet form $(\E^A,D(\E^A))$ is conservative.  
\end{itemize}
\end{thm}
\proof
(i) Let $d=2$. Then by \eqref{eq;intldis} 
\[
\int_{1}^{\infty} \frac{r}{ |\tilde{B}_{r}(0)| } \, dr \ge
\int_{1}^{\infty} \frac{r}{|B_{c_0 r}(0)| } \, dr   =\infty,
\]
where $c_0$ is the constant as in \eqref{eq;intldis}. Therefore by \cite[Theorem 3.4]{St5},  $(\E^A,D(\E^A))$ is recurrent.\\
(ii) Similarly, using \eqref{eq;intldis} one can show that for any $d \ge 2$
\[
\int_{1}^{\infty} \frac{r}{\log \Big (|\tilde{B}_{r}(0)|  \Big )} \, dr \ge
\int_{1}^{\infty} \frac{r}{\log \Big (|B_{c_0 r}(0)|  \Big )} \, dr   =\infty.
\]
Hence by \cite[Theorem 3.6]{St5},  $(\E^A,D(\E^A))$ is conservative.
\qed

By \eqref{eq;intldis}, the completeness  property holds and the doubling property holds since the reference measure is the Lebesgue measure. The weak Poincar\'{e} inequalities on intrinsic balls is also satisfied (see  \cite[Section 1. (b)]{BM}, \cite{Je}, and \cite{JeSa}). Hence the properties (Ia)-(Ic) of \cite{St3} are satisfied. Therefore by \cite[p. 286 A)]{St3} there exists a jointly continuous transition kernel density  $p_{t}(x,y)$ such that
\[
P_t f(x) := \int_{\R^d} p_t(x,y)\, f(y) \ dy, \ \  t>0, \ x,y \in \R^d, \ f\in \mathcal{B}_b(\R^d)
\]
is a $dx$-version of $T_t f$ if $f  \in  L^2(\R^d , dx)_b$. Throughout this paper we set $P_0 : = id$.
Taking the Laplace transform of $p_{\cdot}(x, y)$, we obtain a $\mathcal{B}(\R^d) \times \mathcal{B}(\R^d)$ measurable non-negative resolvent kernel density $r_{\alpha}(x,y)$ such that
\[
R_{\alpha} f(x) := \int_{\R^d} r_{\alpha}(x,y) \ f(y) \ dy,  \ \alpha>0, \ x \in \R^d, \  f \in \mathcal{B}_b(\R^d),
\]
is a $dx$-version of $G_{\alpha} f$ if $f  \in  L^2(\R^d, dx)_b$. For a signed Radon measure $\mu$ on $\R^d$, let us define
\begin{equation*}
R_{\alpha} \mu (x) = \int_{\R^d} r_{\alpha}(x,y) \ \mu(dy), \ \alpha>0, \ x \in \R^d,
\end{equation*}
whenever this makes sense. In particular, $R_{1} \mu$ is a version of $U_{1}\mu$ (see e.g. \cite[Exercise 4.2.2]{FOT}).\\
It follows from \cite[Corollary 4.2 and Remarks (ii) in p.286]{St3} that  for $x, y \in \R^d$, $t > 0$, any $\delta > 0$
\begin{equation}\label{eq;hkestintr}
p_t(x,y) \le c\  \big | \tilde{B}_{\sqrt{t}} (x) \big|^{-1/2}  \big| \tilde{B}_{\sqrt{t}} (y) \big|^{-1/2} \ \exp{ \left ( -\frac{d(x,y)^2}{ (4 + \delta) t} \right)},
\end{equation}
where $c$ is some constant.  

\begin{lemma}\label{l;l1loc}
Let $r > 0$ and $t > 0$. Then 
\[
\sup_{x \in B_r (0)} p_t(x, \cdot) \in L^1 (\R^d,dz).
\]
\end{lemma}
\proof
Let $x,y \in \R^d$ and $t,r > 0$. Note that by \eqref{eq;intudisloc}, $\inf_{x \in \bar{B}_r (0)}  \big|\tilde{B}_{\sqrt{t}}(x)\big| = : M_{t,r} >0$. Putting \eqref{eq;intldis} into \eqref{eq;hkestintr} we obtain
for $x,y \in \R^d$, $t > 0$, any $\delta > 0$
\begin{equation}\label{eq;hkeestno}
p_t(x,y) \le c \ \big|  \tilde{B}_{ \sqrt{t}  }(x) \big|^{-1/2} \  \big| \tilde{B}_{\sqrt{t}}(y) \big|^{-1/2}  \ \exp\left(- \frac{\| x-y \|^{2}}{c_0^2 (4 + \delta) t} \right).
\end{equation}
Using the doubling property, \eqref{eq;hkeestno} can be rewritten as
\[
p_t(x,y) \le c_1\   \frac{1}{ \big| \tilde{B}_{  \sqrt{t}  }(x) \big|} \ \exp{ \left ( -\frac{\| x-y \|^2}{ c_0^2 (4 + \delta) t} \right)},
\]
where $c_1$ is some constant (cf. \cite[proof of Lemma 3.2]{ShTr13a} and \cite[p. 287]{St3}). Therefore
\[
\sup_{x \in B_r (0)} p_t(x, \cdot) \in L^1 (\R^d,dz).
\]
\qed

Using Lemma \ref{l;l1loc}, we show that $(P_t)_{t \ge 0}$ is strong Feller:
\begin{prop}\label{pr;stfellsu}
$(P_t)_{t \ge 0}$ (resp. $(R_{\alpha})_{\alpha > 0}$) is strong Feller, i.e. for $t>0$, $P_t(\mathcal{B}_b(\R^d)) \subset C_b(\R^d)$ (resp. for $\alpha>0$, $R_{\alpha}(\mathcal{B}_b(\R^d)) \subset C_b(\R^d)$).
\end{prop}
\proof
Let $x_n \to x$ in $\R^d$ as $n \to \infty$. For  $f \in \mathcal{B}_b(\R^d)$ and $t > 0$
\[
| P_t f(x_n) - P_t f(x) | \le \int_{\R^d} | p_t(x_n,y) - p_t(x,y) | \ | f(y) | \ dy,
\]   
which converges to $0$ by Lebesgue and Lemma \ref{l;l1loc} and the continuity of $p_t(\cdot,y)$. Note that clearly for $f \in \mathcal{B}_b(\R^d)$, $t > 0$, $P_t f $ is bounded.  Therefore $(P_t)_{t \ge 0}$ is strong Feller. Since $R_{\alpha} f (x) = \int_0^{\infty} e^{-t} \ P_t f(x) \ dt$ and $\|P_t f\|_{L^{\infty}(\R^d, dx)} \le \|f\|_{L^{\infty}(\R^d, dx)}$ for any $f \in \mathcal{B}_b(\R^d)$, $(R_{\alpha})_{\alpha>0}$ is clearly also strong Feller by Lebesgue.
\qed

\begin{remark}\label{fesubk}
We do not know whether the transition function $(P_t)_{t \ge 0}$ is a Feller semigroup or not.
\end{remark}

Now we adopt the construction method of a Hunt process associated with a given Dirichlet form as introduced in \cite[Section 2]{ShTr13a}. In \cite[Section 2]{ShTr13a} we considered a symmetric, strongly local, regular Dirichlet form $(\E,D(\E))$ on $L^{2}(E, \mu)$ with generator $(L,D(L))$ admitting carr\'e du champ, where $E$ is a locally compact separable metric space and $\mu$ is a positive Radon measure on $(E, \mathcal{B}(E))$ with full support on $E$.\\

There, with the corresponding semigroup $(T_t)_{t > 0}$, the transition function $(P_t)_{t > 0}$, the resolvent kernel $R_1$ w.r.t. $(\E,D(\E))$, and so on, we assumed : 
\begin{itemize}
\item[$(\bf{H1})$] There exists a $\mathcal{B}(E) \times \mathcal{B}(E)$ measurable non-negative map $p_{t}(x,y)$ such that
\[
P_t f(x) := \int_{E} p_t(x,y)\, f(y) \, \mu(dy) \,, \; t>0, \ \ x \in E,  \ \ f \in \mathcal{B}_b(E),
\]
is a (temporally homogeneous) sub-Markovian transition function (see \cite[1.2]{CW}) and an $\mu$-version of $T_t f$ if $f  \in  L^2(E , \mu)_b$.
\item[${(\textbf{H2})^{\prime}}$] We can find $\{ u_n \ | \ n \ge 1 \} \subset D(L) \cap C_0(E)$ satisfying:
\begin{itemize}
\item[(i)] For all $\varepsilon \in \Q \cap (0,1)$ and
$y \in D$, where $D$ is any given countable dense set in $E$, there exists $n \in \N$ such that $u_n (z) \ge 1$, for all $z \in \overline{B}_{\frac{\varepsilon}{4}}(y)$ and $u_n \equiv 0$ on $E \setminus B_{\frac{\varepsilon}{2}}(y)$.
\item[(ii)] $R_1\big( [(1 -L) u_n]^+ \big)$, $R_1\big( [(1 -L) u_n]^- \big)$, $R_1 \big( [(1-L_1)u_n^2]^+ \big)$, $R_1 \big( [(1-L_1)u_n^2]^- \big)$ are continuous on $E$ for all $n \ge 1$ where $L_1$ denotes the $L^1(E, \mu)$-generator of $(\E,D(\E))$.
\item[(iii)] $R_1 C_0(E) \subset C(E)$.
\item[(iv)] For any $f \in C_0(E)$ and $x \in E$, the map $t \mapsto P_t f(x)$ is right-continuous on $(0,\infty)$.
\end{itemize}
\end{itemize}

Under  $(\bf{H1})$ and ${(\textbf{H2})^{\prime}}$ we showed that there exists a Hunt process with $(P_t)_{t \ge 0}$ as transition function (see \cite[Lemma 2.9]{ShTr13a}). \\

We intend to do the same here in our concrete situation, i.e. we first show that $(\E^A,D(\E^A))$ satisfies $(\bf{H1})$ and ${(\textbf{H2})^{\prime}}$ and so finally can construct a Hunt process
\begin{equation*}
\bM: = (\Omega , \mathcal{F}, (\mathcal{F}_t)_{t\geq0}, \zeta ,(X_t)_{t\geq0} , (\P_x)_{x \in \R^d_{\Delta} })
\end{equation*}
 satisfying the absolute continuity condition (as stated in \cite[p. 165]{FOT}) with the transition function $(P_t)_{t \ge 0}$. Here $\Delta$ is the cemetery point, $\R^d_{\Delta} : = \R^d \cup \{\Delta\}$ and the lifetime $\zeta : = \inf \{t \ge 0 \ | \  X_t \in \{\Delta\} \}$ and $P_t(x,B) : = P_t 1_B (x) = \P_x(X_t \in B)$ for any $x \in \R^d$, $B \in \mathcal{B}(\R^d)$, $t \ge 0$. Let $(L^A, D(L^A))$ be the generator of $(\E^A, D(\E^A))$.  Since $a_{ij}, \partial_{j} a_{ij} \in L^{\infty}_{loc}(\R^d,dx)$, we have for $f \in C_0^{\infty}(\R^d)$
\begin{equation}\label{eq;dgenerator}
f \in D(L^A) \quad \text{and} \quad L^A f =  \frac{1}{2}\sum_{i,j = 1}^{d} \left(  a_{ij}  \  \partial_{ij} f +  \partial_{j} a_{ij}  \ \partial_i f \right)  \in L^{\infty} (\R^d,dx)_0. 
\end{equation}

\begin{thm}
There exists a Hunt process $\bM$ 
 satisfying the absolute continuity condition  with transition function $(P_t)_{t \ge 0}$.    
\end{thm}
\proof
Using the transition density estimate  \eqref{eq;hkeestno}, we can see as in  \cite[Proposition 3.3 (ii)]{ShTr13a} that $(\textbf{H1})$ and \textbf{(H2)}$^{\prime}$ (iii), (iv)  hold. Clearly we can find $(u_n)_{n \ge 1}  \subset C_0^{\infty}(\R^d) \subset D(L^A)$ such that $(\textbf{H2})^{\prime}$ (i) is satisfied. Furthermore $(\textbf{H2})^{\prime}$ (ii) for $(u_n)_{n \ge 1}$ satisfying $(\textbf{H2})^{\prime}$ (i) follows from \eqref{eq;dgenerator} and Proposition \ref{pr;stfellsu}.
\qed

\begin{remark}\label{r;pconsh}
By Theorem \ref{t;subconsv} and Proposition \ref{pr;stfellsu}, $\P_x(\zeta = \infty) = 1$ for all $x \in \R^d$.
\end{remark}

\begin{lemma}\label{l;subsmom}
Assume (A1) and (A2) hold. For any relatively compact open set $G\subset \R^d$, 
\[
1_{G} \cdot  a_{ii} \ dx \in S_{00}, \quad 1_{G} \cdot   |\partial_j a_{ij} |  \ dx \in S_{00}.
\] 
\end{lemma}
\proof
For any relatively compact open set $G\subset \R^d$ $ 1_{G} \cdot  a_{ii} \ dx$  and $1_{G} \cdot |\partial_j a_{ij} |  \ dx$ are positive finite measures on $\R^d$. Furthermore by (A1) and Proposition \ref{pr;stfellsu}, $R_1 (1_{G} \cdot   a_{ii} \ dx) \in C_b(\R^d)$ and $R_1(1_{G} \cdot   |\partial_j a_{ij} |  \ dx) \in C_b(\R^d)$. Consequently $1_{G} \cdot  a_{ii}  \ dx  \in S_{00}$, $1_{G} \cdot   |\partial_j a_{ij} |  \ dx \in S_{00}$ (see \cite[Proposition 2.12]{ShTr15}).

\qed

We will refer to \cite{FOT} till the end, hence some of its standard notations may be adopted below without definition.
Let  $f^i (x):= x_i$, $i=1,\dots, d$, $x = (x_1,\dots, x_d)\in \R^d$, be the coordinate functions. Then $f^i \in D(\E^A)_{b,loc}$ and for any $g \in C_0^{\infty}(\R^d)$ the following integration by parts formula holds:
\begin{equation}\label{eq;ibpsubel}
- \E^A(f^i,g)=    \frac{1}{2} \int_{\R^d} \Big( \sum_{j=1}^{d}  \partial_j a_{ij}  \Big) g \ dx, \quad 1 \le i \le d.
\end{equation}
The proof of next theorem is basically similar to \cite[Theorem 3.9]{ShTr15}. But we add the proof for the convenience of readers.
\begin{thm}\label{th;exsolsu}
Assume (A1)-(A2) hold. Then it holds $\P_x$-a.s. for any $x \in \R^d$, $i=1,\dots,d$
\begin{equation}\label{sfdsubea}
X_t^i = x_i + \sum_{j=1}^d \int_0^t  \sigma_{ij} (X_s) \ dW_s^j +   \frac{1}{2}\sum_{j=1}^d \int^{t}_{0}   \partial_j a_{ij} (X_s) \, ds, \quad t \ge 0,
\end{equation}
where $\sqrt{A} = (\sigma_{ij})_{1 \le i,j \le d}$ is the square root of the matrix $A$, $W = (W^1,\dots,W^d)$ is a standard d-dimensional Brownian motion on $\R^d$.
\end{thm}
\proof
By Lemma \ref{l;subsmom}, \eqref{eq;ibpsubel}, and \cite[Theorem 5.5.5]{FOT}, the strict continuous additive functional, locally of zero energy and corresponding to the coordinate function $f^i \in D(\E^A)_{b,loc}$, is given by 
\[
N^{[f^i]}_t =  \frac{1}{2}  \int^{t}_{0}  \left( \sum_{j=1}^d  \partial_j a_{ij} \right) (X_s) \, ds, \quad t \ge 0, \quad 1 \le i \le d. 
\] 
The energy measure of $f^i$ denoted by $\mu_{\langle f^i \rangle}$ satisfies $\mu_{\langle f^i \rangle} =  a_{ii} \ dx$.
By Lemma \ref{l;subsmom} for any relatively compact open set $G\subset \R^d$, $1_{G} \cdot \mu_{\langle f^i \rangle} \in S_{00}$ and so the positive continuous additive functional in the strict sense corresponding to the Reuvz measure $\mu_{\langle f^i \rangle}$ is given by 
\[
\langle M^{[f^i]}  \rangle_t = \int_0^t  a_{ii}(X_s) \ ds,
\]
where $M^{[f^i]}_t$ is the continuous local martingale additive functional in the strict sense corresponding to $f^i$. Furthermore since the covariation is
\[
\quad \langle M^{[f^i]}, M^{[f^j]}  \rangle_t = \int_0^t  a_{ij}(X_s) \ ds,
\]
we can construct a d-dimensional Brownian motion $W$ (on a possibly enlarged probability space $(\overline{\Omega}, \overline{\mathcal{F}}, \overline{\P}_x)$, see \cite[Chapter 3, Theorem 4.2]{KS}), that we call again w.l.o.g. $(\Omega,\mathcal{F}, \P_x)$) such that 
\[
M^{[f^i]}_t = \sum_{j=1}^d \int_0^t \sigma_{ij} (X_s) \ dW_s^j,
\]
where $(\sigma_{ij})_{1 \le i,j \le d} =  \sqrt{A} $ is the square root of the matrix $A$.
Note that the equation \eqref{sfdsubea} holds for all $t \ge 0$ because $(\E^A,D(\E^A))$ is conservative (see Remark \ref{r;pconsh}).
\qed

%\begin{remark}
%In (A1) we could not assume more generally $0 \le \langle A(x) \ \xi, \xi \rangle,  \quad \forall \xi \in \R^d$, $m$-a.e. $x \in \R^d$
%because we need SOMETHING the existence of $\sqrt{A} = (\sigma_{ij})_{1 \le i,j \le d}$ as in Theorem \ref{th;exsolsu}.
%\end{remark}

\section{Degenerate elliptic forms with Lebesgue measure}\label{se3deel}

In this section we consider the following assumption:
\begin{itemize}
\item[(A3)]
Let $B:=(b_{ij})_{1 \le i,j \le d}$ be an elliptic symmetric matrix on $\R^d$, i.e. $b_{ij}(x)$ are Borel measurable functions on $\R^d$
and there exist $\lambda_1, \lambda_2 \in C_b(\R^d)$ with  $0 < \lambda_1 \le \lambda_2$ such that for $dx$-a.e. $x \in \R^d$
\begin{equation}\label{eq;coela}
\lambda_1 (x) \ \| \xi \|^2 \le \langle B(x) \xi, \xi \rangle \le \lambda_2(x) \ \|\xi\|^2, \quad \forall \xi \in \R^d.
\end{equation} 
\end{itemize}
We say this matrix $B$ is degenerate (or locally uniformly elliptic) since it can not be uniformly bounded away from zero in \eqref{eq;coela}. Now we fix a matrix $B=(b_{ij})_{1 \le i,j \le d}$ satisfying (A3) and consider the symmetric bilinear form
\[
\E^B(f,g) : = \frac{1}{2} \int_{\R^d} \langle B \ \nabla f , \nabla g \rangle \ dx, \quad f,g \in C_0^{\infty}(\R^d).
\]
By \cite[Chapter II. 2. b)]{MR} $(\E^B,C_0^{\infty}(\R^d))$ is closable in $L^2(\R^d, dx)$ and its closure $(\E^B,D(\E^B))$ is a strongly local, regular, symmetric Dirichlet form. As before in Section \ref{sepo2}, we denote  the $L^2(\R^d, dx)$-semigroup and resolvent associated to $(\E^B,D(\E^B))$ by $(T_t)_{t > 0}$ and $(G_{\alpha})_{\alpha > 0}$. Correspondingly, we can define the intrinsic metric $d(\cdot,\cdot)$ and the intrinsic balls $\tilde{B}_r(x)$,  $x \in \R^d$, $r>0$ relevant to $(\E^B,D(\E^B))$ as introduced in Section \ref{sepo2}.
Furthermore, we assume the (scaled) weak Poincar\'e inequality:
\begin{itemize}
\item[(A4)]
There exists a constant $c>0$ such that 
\[
\int_{\tilde{B}_{r} (x)} | u - \tilde{u}_{x,r} |^2  \ dy  \le  c r^2 \int_{\tilde{B}_{2r}(x)} \langle B \ \nabla u , \nabla u \rangle \ dy, \quad \forall u \in C^{\infty}_0(\R^d), x \in \R^d, r>0,
 \]
where $\tilde{u}_{x,r} = \frac{1}{| \tilde{B}_r(x) |}\int_{\tilde{B}_r(x)} u \ dy$.
\end{itemize}

\begin{remark}\label{r;h12h5}
 Suppose that the symmetric matrix $B = (b_{ij})_{1 \le i,j \le d}$ satisfies (A1) and (A2). Then this $B$ satisfies (A4) (see Section \ref{sepo2}).
\end{remark}

\begin{lemma}\label{le;l;intrinsicdell}
For any $x,y \in \R^d$
\begin{equation}\label{eq;intllowpo}
d (x,y) \ge \frac{1}{\sqrt{c_2}} \ \| x-y \|, \quad  c_2 : =  \sup_{x \in \R^d} \lambda_2(x),
\end{equation}
and  for any bounded set $D \in \mathcal{B}(\R^d)$
\begin{equation}\label{eq;intudisupo}
d(x,y)  \le \frac{1}{\sqrt{c_D}} \| x - y \|, \quad x,y \in D,
\end{equation} 
where $c_D:= \inf_{x \in D} \lambda_1(x)$.
\end{lemma}
\proof
We basically follow the ideas in the proof of \cite[Theorem 4.1]{St5}. For any $z \in \R^d$ the map
\[
u \ : \ x \longmapsto u(x) : = \langle x, z \rangle 
\]
lies in $D({\E^B})_{loc} \cap C(\R^d)$. For fixed $y, y^{\prime} \in \R^d$, $y \neq y^{\prime}$, choose 
\[
z = \frac{(y - y^{\prime})}{\sqrt{c_2} \ \| y - \ y^{\prime} \|} \in \R^d, \quad  c_2 : =  \sup_{x \in \R^d} \lambda_2(x).
\] 
Then by \eqref{eq;coela}
\[
\int_A  d \Gamma^{B}(u,u)  = \int_A  \langle B \nabla u,  \nabla u \rangle dx \le c_2 \int_A  \| \nabla u\|^2  \ dx= \int_A   dx, \quad \forall A \in \mathcal{B}(\R^d).
\]
Hence $\Gamma^{B}(u,u)   \le  dx$.
Furthermore
\[
u(y) - u(y^{\prime}) = \frac{1}{\sqrt{c_2}} \ \| y- y^{\prime} \|.
\]
Therefore for any $x,y \in \R^d$
\[
d (x,y) \ge \frac{1}{\sqrt{c_2}} \ \| x-y \|.
\]
Conversely, let $u \in D(\E^B)_{loc} \cap C(\R^d)$ with
\begin{equation}\label{eq;absosub}
\Gamma^{B}(u,u) \le dx.
\end{equation}
Let $( u_n )_{n \ge 1} \subset C_0^{\infty}(\R^d)$ be a sequence that converges to $u$ locally in $\sqrt{\E_1^B}$-norm and locally uniformly. Then by  \eqref{eq;coela}, $\left(\partial_i u_n   \right)_{n \ge 1}$ is a Cauchy sequence in $L^2_{loc}(\R^d, dx)$. Therefore there exists $v_i \in L^2_{loc} (\R^d, dx)$, $i= 1, \dots, d$ such that $\partial_i u_n  \to v_i$ in $L^2_{loc}(\R^d, dx)$ as $n \to \infty$. 
Let $D \in \mathcal{B}(\R^d)$ be a bounded set. Then by \eqref{eq;coela}
\[
\lim_{n \to \infty} \int_{D} \langle B \nabla u_n,    \nabla u_n \rangle  \ dx = \int_{D} \langle B v, v \rangle \ dx.
\]
Then
\[
\int_{D} \ d \Gamma^B (u,u) = \lim_{n \to \infty} \int_{D} \ d \Gamma^B (u_n,u_n)  = \lim_{n \to \infty} \int_{D} \langle B \nabla u_n,    \nabla u_n \rangle  \ dx = \int_{D} \langle B v, v \rangle \ dx,
\]
where $v= (v_1, \dots, v_d).$ Together with \eqref{eq;absosub} this implies that 
\[
\int_{D} \langle B v, v \rangle \ dx  \le  \int_{D} 1  \ dx.
\] 
In particular, by \eqref{eq;coela}
\[
c_D \ \| v \|^2 \le 1
\]
$dx$-a.e. on $D$ where $c_D:= \inf_{x \in D} \lambda_1(x)$. Now  following the proof of \cite[Theorem 4.1, p.264]{St5} one can show that
\[ 
|u(x) - u(y)| \le \frac{1}{\sqrt{c_D}} \| x - y \|,  \quad \forall x,y \in D.
\]
\qed

\begin{remark}
By Lemma \ref{le;l;intrinsicdell}  the topology induced by the intrinsic metric coincides with the Euclidean topology.  Hence the Dirichlet form $(\E^B,D(\E^B))$ is strongly regular.
\end{remark}

By \eqref{eq;intllowpo}, the completeness  property holds.  The doubling property holds since the reference measure is the Lebesgue measure. By the assumption (A4) the weak Poincar\'{e} inequality on intrinsic balls is also satisfied. Hence the properties (Ia)-(Ic) of \cite{St3} are satisfied. Therefore likewise Section \ref{sepo2}  by \cite[p. 286 A)]{St3} there exists a jointly continuous transition kernel density  $p_{t}(x,y)$ such that
\[
P_t f(x) := \int_{\R^d} p_t(x,y)\, f(y) \ dy, \ \  t>0, \ x,y \in \R^d, \ f\in \mathcal{B}_b(\R^d)
\]
is a $dx$-version of $T_t f$ if $f  \in  L^2(\R^d , dx)_b$. Furthermore,
it follows from \cite[Corollary 4.2 and Remarks (ii) in p.286]{St3} that  for $x, y \in \R^d$, $t > 0$, any $\delta>0$
\begin{equation}\label{eq;hkestpo}
p_t(x,y) \le c\  \big| \tilde{B}_{\sqrt{t}} (x) \big|^{-1/2}  \big| \tilde{B}_{\sqrt{t}} (y) \big|^{-1/2} \ \exp{ \left ( -\frac{d(x,y)^2}{ (4 + \delta) t} \right)},
\end{equation}
where $c$ is some constant.  Similarly, $R_{\alpha} f$ and $R_{\alpha} \mu$ can be defined as in Section \ref{sepo2}.
\begin{thm}\label{t;pocovat}
\begin{itemize}
\item[(i)] Let $d=2$. Then the Dirichlet form $(\E^B,D(\E^B))$ is recurrent.
\item[(ii)] The Dirichlet form $(\E^B,D(\E^B))$ is conservative.  
\end{itemize}
 \end{thm}
\proof
Using \eqref{eq;intllowpo} the proof is similar to Theorem \ref{t;subconsv}.
\qed

\begin{lemma}\label{l;l1locell}
Let $t,r > 0$. Then 
\[
\sup_{x \in B_r (0)} p_t(x, \cdot) \in L^1 (\R^d,dz).
\]
\end{lemma}
\proof
Using \eqref{eq;intllowpo}, \eqref{eq;hkestpo}, and the doubling property, the proof is similar to Lemma \ref{l;l1loc}.
\qed

\begin{prop}\label{pr;stfeullpo}
$(P_t)_{t \ge 0}$ and $(R_{\alpha})_{\alpha > 0}$ are strong Feller (cf. Proposition \ref{pr;stfellsu}).
\end{prop}
\proof
Using Lemma \ref{l;l1locell} the proof is similar to Proposition \ref{pr;stfellsu}. So we omit it.
\qed

In contrast to the case of the subelliptic Dirichlet form $(\E^A,D(\E^A))$ considered in Section \ref{sepo2} we obtain: 
\begin{thm}\label{t;fesepo}
The transition function $(P_t)_{t \ge 0}$ satisfies:
\begin{itemize}
\item[(i)] $\lim_{t \to 0} P_t f(x) = f(x)$ for each $x \in \R^d$ and $f \in C_0(\R^d)$.
\item[(ii)] $P_t C_0(\R^d) \subset C_{\infty}(\R^d)$ for each $t>0$.
\end{itemize} 
In particular, $(P_t)_{t \ge 0}$ is a Feller semigroup.  
\end{thm}
\proof
By \eqref{eq;intudisupo} there exists a constant $c_x > 0$ (depending on $t$ and $x$) such that 
\[
| B_{ (c_x \sqrt{t}) } (x) | \le    | \tilde{B}_{\sqrt{t}} (x) | 
\]
and a constant $c_y > 0$ (depending on $t$ and $y$) such that 
\[
|B_{ (c_y \sqrt{t}) } (y) |  \le   | \tilde{B}_{\sqrt{t}} (y) |. 
\]
Therefore together with the doubling property \eqref{eq;hkestpo} can be rewritten as  
\begin{equation}\label{eq;heust}
p_t(x,y) \le c_1 \ \frac{1}{| B_{ (c_x \sqrt{t}) } (x) |}  \ \exp\left(- \frac{\| x-y \|^{2}}{c_2^2 (4 + \delta) t} \right),
\end{equation}
and using symmetry of $p_t(\cdot,\cdot)$
\[
p_t(x,y) \le c_1 \ \frac{1}{|B_{ (c_y \sqrt{t}) } (y) |}  \ \exp\left(- \frac{\| x-y \|^{2}}{c_2^2 (4 + \delta) t} \right),
\]
where $c_1$ is some constant and $c_2$ is the constant as in \eqref{eq;intllowpo}. Note that since $(\E^B,D(\E^B))$ is conservative and $(P_t)_{t \ge 0}$ is strong Feller, we have $P_t 1(x) = 1$ for all $x \in \R^d$, $t >0$. Then for $f \in C_0(\R^d)$, $x \in \R^d$, $t> 0$
\begin{eqnarray*}
\big | P_t f(x) - f(x) \big | &=& \left| \int_{\R^d} p_t(x,y) \big( f(y) - f(x) \big) \ dy \right | \\ 
&\le& c_1 \ \int_{\R^d} \frac{1}{|B_{ (c_x \sqrt{t}) } (x) |}  \ \exp\left(- \frac{\| x-y \|^{2}}{c_2^2 (4 + \delta)   t} \right) \ \big |  f(y) - f(x) \big | \ dy, 
\end{eqnarray*}
which converges to zero as $t$ tends to zero.
Furthermore for $f \in C_0 (\R^d) $, $x \in \R^d$, $t>0$ by Proposition \ref{pr;stfeullpo} 
\[
P_t f  \in C(\R^d)
\]
and 
\[
P_t f(x) = \int_{\R^d} p_t(x,y) f(y) \ dy \le c_1 \ \int_{\R^d} \frac{1}{|B_{ (c_y \sqrt{t}) } (y) |}  \ \exp\left(- \frac{\| x-y \|^{2}}{c_2^2 (4 + \delta)  t} \right) \ f(y) \ dy, 
\]
which converges to zero as $\|x\| $ goes to infinity. In particular, by \cite[Lemma 2.3]{ShTr13a} $(P_t)_{t \ge 0}$ is a Feller semigroup.
\qed

\begin{remark}
Under the assumptions of (A1) and (A2) in Section \ref{sepo2} we do not know whether the transition function
$(P_t)_{t \ge 0}$ associated with $(\E^A,D(\E^A))$ in Section \ref{sepo2} is a Feller semigroup or not (cf. Remark \ref{fesubk}). However Theorem \ref{t;fesepo} says that if we add the assumption (A3), the transition function $(P_t)_{t \ge 0}$ is a Feller semigroup (cf. Remark \ref{r;h12h5}).  
\end{remark}

According to  Theorem \ref{t;fesepo} and the classical Feller theory, there exists a Hunt process 
\[
\bM = (\Omega , \mathcal{F}, (\mathcal{F}_t)_{t\geq0}, \zeta ,(X_t)_{t\geq0} , (\P_x)_{x \in \R^d_{\Delta} })
\]
satisfying the absolute continuity condition with the transition function $(P_t)_{t \ge 0}$.

\begin{remark}\label{r;pconshepo}
By Theorem \ref{t;pocovat} and Proposition \ref{pr;stfeullpo}, $\P_x(\zeta = \infty) = 1$ for all $x \in \R^d$.
\end{remark}
For the time being, we consider a rather strong assumption:
\begin{itemize}
\item[(A5)] For each $i,j =1 , \dots ,d$, $\partial_j  b_{ij}  \in L_{loc}^{\infty}(\R^d,dx)$.
\end{itemize}

\begin{lemma}
Assume (A3)-(A5). Then for any relatively compact open set $G\subset \R^d$, 
\[
1_{G} \cdot  b_{ii}  \ dx  \in S_{00}, \quad 1_{G} \cdot   |\partial_j b_{ij} |  \ dx \in S_{00}.
\] 
\end{lemma}
\proof
For any relatively compact open set $G\subset \R^d$ $ 1_{G} \cdot  b_{ii} dx$  and $1_{G} \cdot |\partial_j b_{ij} |  dx$ are positive finite measures on $\R^d$. Furthermore by (A3), (A5), and Proposition \ref{pr;stfeullpo}, $R_1 (1_{G} \cdot   b_{ii} dx) \in C_b(\R^d)$ and $R_1(1_{G} \cdot   |\partial_j b_{ij} |  dx) \in C_b(\R^d)$. Therefore $1_{G} \cdot  b_{ii}  dx  \in S_{00}$ and $1_{G} \cdot   |\partial_j b_{ij} |  dx \in S_{00}$ (see \cite[Proposition 2.12]{ShTr15}). 
\qed

\begin{remark}
It is not possible to obtain the resolvent density estimate by taking Laplace transform in \eqref{eq;heust} w.r.t $t$ because the constant $c_x$ in \eqref{eq;heust} depends on $t$. This is the reason why we assume local boundedness of $\partial_j  b_{ij}$ as in (A5). However, the assumption (A5) can be relaxed in the next subsection.
\end{remark}

Finally, we obtain:
\begin{thm}\label{t;sdpoe}
Assume (A3)-(A5). It holds $\P_x$-a.s. for any $x \in \R^d$, $i=1,\dots,d$
\begin{equation}\label{eq;fupoijd}
X_t^i = x_i + \sum_{j=1}^d \int_0^t  \rho_{ij} (X_s) \ dW_s^j +   \frac{1}{2}\sum_{j=1}^d \int^{t}_{0}   \partial_j b_{ij} (X_s) \, ds, \quad t \ge 0,
\end{equation}
where $\sqrt{B} = (\rho_{ij})_{1 \le i,j \le d}$ is the positive square root of the matrix $B$, $W = (W^1,\dots,W^d)$ is a standard d-dimensional Brownian motion on $\R^d$.
\end{thm}
\proof
The proof is simialr to Theorem \ref{th;exsolsu}. So we omit it.
\qed

\subsection{Nash-type inequality and part processes}
This subsection is devoted to relaxing the strong assumption (A5) to a rather weak assumption (A5)$^{\prime}$ below. We mainly use the Nash-type inequality and part Dirichlet forms of $(\E^B,D(\E^B))$. For the notations (especially concerning part forms and part processes) which appear in this subsection we refer to \cite[Section 2]{ShTr13a} (cf. \cite{FOT}).  Let 
\begin{equation*}
B_k : = \{x \in \R^d  \ | \ \|x \| < k \}, \quad k \ge 1,
\end{equation*}
and for any $G \subset \R^d$
\[
C^{\infty}(\overline{G}) : = \{f:\overline{G} \to \R \ | \ \exists g \in C_0^{\infty}(\R^d), g|_{\overline{G}} = f   \}.
\]
According to \eqref{eq;coela} the closure of 
\[
\E^{B,\overline{B}_k}(f,g) := \frac{1}{2} \int_{B_k}  \langle  B \nabla f, \nabla g \rangle \ dx, \quad f,g \in C^{\infty}(\overline{B}_k),
\]
in $L^2(\overline{B}_k,dx) \equiv L^2(B_k,dx)$, $k \ge 1$, denoted by $(\E^{B,\overline{B}_k},D(\E^{B,\overline{B}_k}))$, is a regular Dirichlet form on $\overline{B}_k$. 
\begin{lemma}\label{l;ssashhke}
\begin{itemize}
\item[(i)] The following Nash-type inequality holds:
\begin{itemize}
\item[(a)] if $d \ge 3$, then for $f \in D(\E^{B,\overline{B}_k})$
\begin{equation*}
\left\|f\right\|_{2,B_k}^{2 + \frac{4}{d}}\leq c_k \left[\E^{B,\overline{B}_k}(f,f) + \left\|f\right\|_{2,B_k}^2 \right]\left\|f\right\|_{1,B_k}^{\frac{4}{d}},
\end{equation*}
\item[(b)] if $d=2$, then for $f \in D(\E^{B,\overline{B}_k})$ and any $\delta>0$
\begin{equation*}
\left\|f\right\|_{2,B_k}^{2 + \frac{4}{d+\delta}} \le c_k \left[\E^{B,\overline{B}_k}(f,f) + \left\|f\right\|_{2,B_k}^2 \right]\left\|f\right\|_{1,B_k}^{\frac{4}{d+\delta}}.
\end{equation*}
Here $c_k >0 $ is a constant which goes to infinity as $k \rightarrow \infty$.
\end{itemize}

\item[(ii)]
We obtain for $m$-a.e. $x, y \in B_k$
\begin{itemize}
\item[(a)] if $d \ge 3$, then
\[
r^{B_k}_{1} (x,y)  \le c_1 \frac{1}{\|x-y\|^{d-2}},
\]
\item[(b)] if $d=2$, then for any $\delta>0$
\[
r^{B_k}_{1} (x,y)  \le c_2 \frac{1}{\|x-y\|^{d+\delta-2}}.
\]
Here $c_1, c_2>0$ are some constants.
\end{itemize}

\end{itemize}
\end{lemma}
\proof
We can apply Sobolev's inequality on each $B_k$. Using  \eqref{eq;coela} we derive the Nash type inequalities in (i) (see \cite[Lemma 5.4]{ShTr13a}). Following  the proof of \cite[Proposition 5.5, Corollary 5.6]{ShTr13a} the assertion (ii) follows.
\qed

Now we replace (A5) by
\begin{itemize}
\item[(A5)$^{\prime}$]   $\partial_j  b_{ij}  \in L_{loc}^{\frac{d}{2} + \varepsilon}(\R^d,dx)$ for some $\varepsilon>0$ and  each $i,j =1 , \dots ,d$.
\end{itemize}
\begin{lemma}\label{l;ssmooloc}
Assume (A3), (A4), and (A5)$^{\prime}$. Let  $f \in L^{\frac{d}{2} + \varepsilon} ( B_k, dx)$ for some $\varepsilon > 0$. Then
\[
1_{B_k} \cdot |f| dx \in S^{B_k}_{00}.
\]
In particular 
\[
1_{B_k} \cdot b_{ii} dx  \in S_{00}^{B_k}, \quad 1_{B_k} \cdot   |\partial_j b_{ij} |  dx \in S_{00}^{B_k}.
\] 
\end{lemma}
\proof
Using the  estimate of resolvent density as in Lemma \ref{l;ssashhke} (ii) and (A5)$^{\prime}$, the proof is similar to the proof of \cite[Lemma 5.8]{ShTr13a}). So we omit it.
\qed

The following integration by parts formula holds for the coordinate functions $f^i\in D(\E^{B,B_k})_{b,loc}$, $i=1,\dots,d$ and $g \in C_0^{\infty}(B_k)$: 
\begin{equation}\label{eq;spibp}
- \E^{B,B_k}(f^i,g)=    \frac{1}{2}  \int_{B_k}  \left(  \sum_{j=1}^{d} \partial_j b_{ij} \right)  g \ dx.
\end{equation}
Let $D_{B_k^c} := \inf\{t \ge 0  \ | \ X_t \in B_k^c\}$.
\begin{prop}\label{t;c5lsfd3}
Assume (A3), (A4), and (A5)$^{\prime}$. Then the process $\bM$ satisfies 
\begin{equation}\label{eq;fspp}
X_t^i = x_i + \sum_{j=1}^d \int_0^t  \rho_{ij} (X_s) \ dW_s^j +   \frac{1}{2}\sum_{j=1}^d \int^{t}_{0}   \partial_j b_{ij} (X_s) \, ds, \quad  t < D_{B_k^c},
\end{equation}
$\P_x$-a.s. for any $x \in B_k$, $i=1,\dots,d$ where $W$ is a standard d-dimensional Brownian motion on $\R^d$.
\end{prop}
\proof
Applying \cite[Theorem 5.5.5]{FOT} to $(\E^{B,B_k}, D(\E^{B,B_k}))$, the assertion then follows from Lemma \ref{l;ssmooloc} and \eqref{eq;spibp} (see Theorem \ref{th;exsolsu} for details).
\qed

\begin{lemma}\label{l;licapli}
 For all $x \in  \R^d$
\[
\P_x \Big(\lim_{k \rightarrow \infty} D_{B_k^c} = \infty \Big)=1.
\]
\end{lemma}
\proof
The proof is similar to \cite[Lemma 5.10]{ShTr13a}. So we omit it. 
\qed

\begin{thm}\label{ch5;t;solex0}
Assume (A3), (A4), and (A5)$^{\prime}$. Then the process $\bM$ satisfies \eqref{eq;fupoijd} for all $x \in \R^d$.
\end{thm}
\proof
Let $k \to \infty$ in \eqref{eq;fspp}. Then by Lemma \ref{l;licapli} the result follows. 
\qed

\begin{remark}
The strict decomposition associated to the Dirichlet form with the uniformly elliptic matrix is presented in \cite[Example]{fuku93}.  Note that the  uniformly elliptic matrix  clearly satisfies (A3) and (A4). Furthermore the assumptions (A5)$^{\prime}$ is weaker than the assumption 
that  $\partial_j  b_{ij} $ is locally bounded as in  \cite[Example]{fuku93}. Therefore the Dirichlet form $(\E^B,D(\E^B))$ includes the case in \cite[Example]{fuku93}.
\end{remark}

\section{Pathwise unique and strong solutions}\label{sect4pust}
In this section we show that the weak solutions appearing in Section \ref{sepo2}, \ref{se3deel} can be pathwise unique and strong solutions.

\begin{prop}\label{ch4;t;ssoleae1}
Assume that (A1) and (A2) hold. Then the (weak) solution in Theorem \ref{th;exsolsu} is strong and pathwise unique. In particular, it is adapted to the filtration $(\mathcal{F}_t^W)_{t\geq0}$ generated by the Brownian motion $(W_t)_{t\geq0}$ as in \eqref{sfdsubea} and its lifetime is infinite.
\end{prop}
\proof
Since $(a_{ij})_{1 \le i,j \le d}$ are twice continuously differentiable, the square root $(\sigma)_{1 \le i,j \le d}$  and  the first derivatives of $(a_{ij})_{1 \le i,j \le d}$ are locally Lipschitz continuous on $\R^d$ (see \cite[Proposition 6.2 (ii) in Chapter 4]{IW}).  
Hence by \cite[Theorem 3.1]{IW} the SDE  \eqref{sfdsubea} has a pathwise unique strong solution up to the explosion time. 
Therefore the (weak) solution in Theorem \ref{th;exsolsu} is strong, pathwise unique and its lifetime  is infinite by Remark \ref{r;pconsh}.
\qed

We additionally assume (in Section \ref{se3deel})
\begin{itemize}
\item[(A6)] for each $1 \le i,j \le d$,
\begin{itemize}
\item[(i)] The matrix $B$  is given by  $B = \rho \rho^{\prime}$ where  $\rho^{\prime}$ is a transpose of the matrix $\rho = (\rho_{ij})_{1 \le i,j \le d}$.
\item[(ii)]   $\rho_{ij}$ is bounded continuous on $\R^d$,
\item[(iii)]   There exists a  constant $c_B>0$  such that $c_B \| \xi \|^2  \le \langle  \rho \rho^{\prime} (x) \ \xi, \xi \rangle$ for all  $x, \xi \in \R^d$,
\item[(iv)]  $ \left \| \nabla  \rho_{ij}  \right \| \in L^{2(d+1)}_{loc} (\R^d,dx)$
\end{itemize}
\end{itemize}

\begin{remark}
\begin{itemize}
\item[(i)] The assumptions  (A6) (ii) and (iv) implies that $\partial_j  b_{ij}  \in L_{loc}^{ 2 (d+1)}(\R^d,dx)$.
\item[(ii)] The assumption (A6) implies (A3), the weak Poincar\'e inequality (A4), and (A5)$^{\prime}$.
\end{itemize}
\end{remark}

\begin{prop}\label{ch4;t;ssoleae2}
Assume that  (A6) holds. Then the (weak) solution in Theorem \ref{t;sdpoe}  is strong and pathwise unique. In particular, it is adapted to the filtration $(\mathcal{F}_t^W)_{t\geq0}$ generated by the Brownian motion $(W_t)_{t\geq0}$ as in \eqref{eq;fupoijd} and its lifetime is infinite.
\end{prop}
\proof
It follows from \cite[Theorem 1.1]{Zh} that for given Brownian motion $(W_t)_{t\geq0}$, $x \in \R^d$ as in \eqref{eq;fupoijd} there exists a pathwise unique strong solution to  \eqref{eq;fupoijd} up to its explosion time. Therefore the (weak) solution in Theorem \ref{t;sdpoe} is strong, pathwise unique and its lifetime  is infinite by Remark \ref{r;pconshepo}.
\qed

%\begin{remark}
%For unique strong solutions to  the SDE \eqref{sdedd}  up to lifetime, \cite[Theorem 1.1]{Zh} presents two non-explosion conditions. The two non-explosion conditions are not satisfied by the assumptions (A1)  and (A6). On the other hand, with the assumptions (A1) and (A6),   by Theorem \ref{ch4;t;ssoleae} and its proof,  we know that the solution to \eqref{sfdsubea} up to its lifetime fits into the frame of \cite[Theorem 1.1]{Zh} and so is unique strong. Furthermore its lifetime is infinite. Therefore, the redundant condition to show strong unique solution of \eqref{sfdsubea} , i.e.
%\[
%a_{ij} \in C^{\infty} (\R^d) \cap C_b(\R^d)  \  \text{in} \ (A1)
%\]
%provides another non-explosion criterion in \cite[Theorem 1.1]{Zh}.
%\end{remark}

\vspace*{2cm}
Jiyong Shin\\
School of Mathematics \\
Korea Institute for Advanced Study\\
85 Hoegiro Dongdaemun-gu,\\
Seoul 02445, South Korea,  \\
E-mail: yonshin2@kias.re.kr \\

\end{document}